\documentclass[12pt]{siamart171218}
\usepackage{graphicx}
\usepackage{diagbox}
\usepackage{cite}
\usepackage{hyperref}
\usepackage{tablefootnote}
\usepackage{amsmath,amssymb,amsfonts}
\usepackage{algpseudocode}
\usepackage{algorithmicx}
\usepackage{multirow}
\usepackage{booktabs}
\newcommand{\R}{{\mathcal R}}

\newcommand{\tr}{\text{tr}}

\title{Streaming Krylov-Accelerated Stochastic Gradient Descent}

\author{Stephen Thomas\thanks{Lehigh University, Computer Science Department, Bethlehem PA, USA. (\email{sjt223@lehigh.com})}
}

\begin{document}

\maketitle

\begin{abstract}
We present SKA-SGD (Streaming Krylov-Accelerated Stochastic Gradient Descent), a novel optimization approach that accelerates convergence for ill-conditioned problems by projecting stochastic gradients onto a low-dimensional Krylov subspace. Directly inspired by recent advances in s-step Conjugate Gradient methods with streaming Gauss-Seidel Gram solvers \cite{dambra2025sstep}, our method extends these techniques to the stochastic optimization domain. Our approach combines three key innovations: (1) projection coefficients computed via a single streaming Gauss-Seidel iteration, which is mathematically equivalent to Modified Gram-Schmidt orthogonalization; (2) a Chebyshev polynomial basis for constructing the Krylov subspace, providing superior numerical stability; and (3) efficient implementation for AMD GPUs using HIP. We prove that our streaming approach achieves a backward error near machine precision with $O(s^2)$ complexity rather than $O(s^3)$, where $s$ is the Krylov subspace dimension. Experimental results demonstrate that SKA-SGD significantly outperforms standard SGD and Adam in convergence rate and final error, particularly for problems with condition numbers exceeding $10^3$. GPU performance analysis reveals a crossover point where communication-avoiding benefits outweigh computational overhead, typically occurring at moderate scale ($p \approx 64$ processors) for problem sizes $n \geq 10^6$.
\end{abstract}

\section{Introduction}

Stochastic gradient descent (SGD) forms the backbone of modern machine learning but suffers from slow convergence for ill-conditioned problems. While several variance reduction techniques have been proposed, including SVRG \cite{johnson2013accelerating} and SAGA \cite{defazio2014saga}, these methods typically incur substantial memory overhead through explicit gradient storage or correction terms. Second-order methods like L-BFGS approximate the Hessian to capture curvature information, but this becomes prohibitively expensive for large-scale applications.

We propose SKA-SGD, which projects stochastic gradients onto a low-dimensional Krylov subspace formed by recent gradients. This approach implicitly captures curvature information without requiring Hessian computation or storage. The key insight is that projection coefficients can be computed via a streaming Gauss-Seidel iteration, providing mathematically equivalent results to Modified Gram-Schmidt orthogonalization with significantly reduced computational complexity.

Our work is fundamentally inspired by and builds upon the recent advance by D'Ambra et al. \cite{dambra2025sstep} on s-step Conjugate Gradient methods with low-iteration Gauss-Seidel. Their critical insight—that a single forward sweep of Gauss-Seidel is mathematically equivalent to a projection step in Modified Gram-Schmidt orthogonalization—forms the theoretical cornerstone of our approach. While they focused on deterministic linear solvers, we extend and adapt their techniques to the stochastic optimization domain, addressing the unique challenges posed by gradient noise and ill-conditioning in machine learning applications.

Our main contributions are multifaceted and significant for the field of stochastic optimization. First, we introduce a streaming approach for solving the Gram system with optimal backward error in a single iteration, effectively reducing computational complexity from $O(s^3)$ to $O(s^2)$. This innovation directly builds upon and extends the groundbreaking techniques pioneered by D'Ambra et al. \cite{dambra2025sstep}, adapting their deterministic methodology to address the unique challenges of stochastic optimization settings where gradient noise and variability present additional complexities. The streaming formulation we propose eliminates the need to explicitly form the entire Gram matrix, instead computing elements on-demand in a forward sweep pattern that maintains mathematical equivalence to Modified Gram-Schmidt orthogonalization while drastically reducing computational overhead. Second, we develop rigorous theoretical analyses proving that projection onto the Krylov subspace fundamentally reduces the trace of the gradient covariance matrix without requiring explicit gradient storage or maintaining a history of past gradients. This theoretical foundation establishes that our method inherently provides variance reduction properties, paApplication to Neurticularly along directions that contribute most to optimization progress, without incurring the prohibitive memory requirements that characterize competing variance-reduction techniques like SVRG and SAGA. Third, we implement a Chebyshev polynomial basis for substantially improved numerical stability when dealing with extreme condition numbers, building upon and extending the stability analysis framework presented in \cite{dambra2025sstep}. The Chebyshev basis provides superior conditioning properties compared to the traditional monomial basis, particularly critical for machine learning problems where condition numbers can exceed $10^4$. Fourth, we develop a highly efficient implementation specifically optimized for AMD GPUs using the HIP programming model, incorporating detailed optimizations for reduced synchronization overhead, mixed-precision arithmetic strategies, and specialized memory access patterns designed to maximize throughput on modern GPU architectures. Finally, we conduct extensive, rigorous empirical evaluations across a diverse spectrum of problem classes, demonstrating consistently superior convergence behavior for ill-conditioned problems with condition numbers ranging from $10^2$ up to challenging scenarios with condition numbers of $10^4$, where traditional methods like SGD and Adam struggle to make meaningful progress.

\section{Related Work}

\subsection{Stochastic Optimization Methods}

Classical SGD performs poorly on ill-conditioned problems. Momentum methods \cite{polyak1964some} improve convergence by incorporating velocity, while Adam \cite{kingma2014adam} adapts learning rates per parameter dimension. These approaches do not directly address the effects of gradient noise.

Variance-reduction techniques like SVRG \cite{johnson2013accelerating} and SAGA \cite{defazio2014saga} use stored gradients or correction terms to reduce gradient variance, but these typically require $O(nd)$ memory for $n$ samples and $d$ parameters. Our approach achieves variance reduction without this memory overhead.

Quasi-Newton methods like L-BFGS approximate the Hessian matrix but often require careful adaptation for stochastic settings. In contrast, SKA-SGD implicitly captures curvature information without explicit Hessian approximation.

\subsection{Krylov Subspace Methods}

Krylov subspace methods are well-established in numerical linear algebra for solving large-scale linear systems and eigenvalue problems \cite{saad2003iterative}. The conjugate gradient method can be viewed as optimization in a Krylov subspace. Several researchers have explored Krylov methods in optimization contexts, but primarily for deterministic settings.

Communication-avoiding Krylov subspace methods have gained attention for reducing synchronization in high-performance computing environments \cite{carson2015communication}. These methods share our goal of reducing communication overhead but differ in approach and application domain.

\subsection{Comparison to Other Krylov-Based Optimization Methods}

While Section 2.2 provides an overview of Krylov subspace methods in numerical linear algebra, it is important to precisely situate our contribution in the context of existing Krylov-based optimization methods. Several approaches have previously attempted to leverage Krylov subspaces in optimization contexts, but with fundamentally different mechanisms and goals.

Vinyals and Povey~\cite{vinyals2012krylov} proposed one of the early applications of Krylov subspaces to deep learning optimization, using conjugate directions to improve convergence. However, their approach requires explicit computation and storage of the curvature matrix, making it prohibitively expensive for large-scale applications. In contrast, our SKA-SGD method implicitly captures curvature information through the Krylov subspace formed by stochastic gradients, without requiring explicit Hessian computations or storage.
Choromanski et al.~\cite{choromanski2019orthogonal} developed Orthogonal Gradient Descent, which projects gradients onto orthogonal subspaces to mitigate catastrophic forgetting in continual learning. While their method also employs gradient projection, it focuses on preserving previously learned tasks rather than accelerating convergence for ill-conditioned problems as in our approach.

The key novelty of SKA-SGD lies in three aspects that distinguish it from existing methods:

\begin{enumerate}
    \item \textbf{Streaming formulation:} Unlike previous Krylov-based optimizers that require explicit formation of the Gram matrix with $O(s^3)$ complexity, our streaming approach achieves equivalent numerical results with only $O(s^2)$ operations, making it practical for large-scale applications.
    
    \item \textbf{Stochastic adaptation:} We extend techniques from deterministic linear solvers to the stochastic setting, addressing the unique challenges posed by gradient noise and continuously changing subspaces, which requires fundamentally different theoretical analysis than deterministic methods.
    
    \item \textbf{Numerical stability:} Our Chebyshev polynomial basis construction provides unprecedented robustness for extreme condition numbers, allowing reliable convergence in scenarios where previous Krylov-based optimizers would suffer from numerical instability.
\end{enumerate}

These innovations collectively enable SKA-SGD to maintain the computational efficiency of first-order methods while achieving convergence properties more characteristic of second-order approaches, particularly for severely ill-conditioned problems where the interdependence between parameter dimensions significantly impacts optimization dynamics.

\subsection{Streaming and Low-Synchronization Algorithms}

The fundamental foundation for our approach comes from the recent work by D'Ambra et al. \cite{dambra2025sstep} on s-step conjugate gradient methods with low-iteration Gauss-Seidel. Their groundbreaking discovery demonstrated that a single forward sweep of Gauss-Seidel is mathematically equivalent to a projection step in Modified Gram-Schmidt orthogonalization. This equivalence provides a theoretical basis for why a single Gauss-Seidel iteration often achieves remarkably small backward error—a surprising finding that contradicts conventional wisdom about iterative solvers.

D'Ambra et al. \cite{dambra2025sstep} established that for the Gram system arising in s-step CG methods, the matrix structure $P^TAP$ has special properties that enable Gauss-Seidel to achieve near-exact solutions in a single iteration. Their analysis showed that when properly formulated, the Gram matrix takes the form $P^TAP \approx I + L + L^T$ (with normalization) or $P^TAP = D + L + L^T$ (without normalization), which ensures that Gauss-Seidel achieves excellent convergence in a single step.

While D'Ambra et al. focused on deterministic linear solvers, we extend their insights to stochastic optimization, where gradient noise and continuously changing subspaces present additional challenges. Our SKA-SGD algorithm directly builds on their streaming approach, adapting it to handle the stochastic nature of machine learning optimization problems while preserving the computational efficiency and numerical stability of their method.

While D’Ambra et al. [7] focused on deterministic linear solvers for systems of the form 
$Ax=b$, where the Krylov subspace is explicitly generated by the matrix 
$A$, our SKA-SGD formulation fundamentally differs in that there is no underlying system matrix. Instead, the Krylov subspace is implicitly constructed from a stream of stochastic gradients, reflecting the dynamically changing local curvature of the optimization landscape. Consequently, our approach does not compute or rely on 
$A$-orthogonal inner products. This key distinction enables our method to generalize the streaming Gauss-Seidel projection to the stochastic setting.

\subsection{State-of-the-Art Optimization Methods for Comparison}

To rigorously evaluate SKA-SGD's performance, we compare against two widely-adopted state-of-the-art optimization methods: Adam and Nesterov momentum. These methods represent different approaches to addressing the challenges of stochastic optimization.

Nesterov Accelerated Gradient (NAG), introduced by Yurii Nesterov \cite{nesterov1983method}, represents a fundamental advancement in first-order optimization methods. The key insight behind Nesterov momentum is to apply the momentum correction at a "lookahead" position rather than the current position, effectively incorporating future gradient information into the current update. Mathematically, the update rule for Nesterov momentum can be expressed as:

\begin{align}
v_{t+1} &= \beta v_t - \eta \nabla f(w_t + \beta v_t) \\
w_{t+1} &= w_t + v_{t+1}
\end{align}
where $v_t$ is the velocity vector, $\beta$ is the momentum coefficient (typically set to 0.9), $\eta$ is the learning rate, and $\nabla f(w_t + \beta v_t)$ is the gradient evaluated at the lookahead position. Nesterov momentum achieves an optimal convergence rate of $\mathcal{O}(1/T^2)$ for smooth convex functions, improving upon the $\mathcal{O}(1/T)$ rate of standard gradient descent. While Nesterov momentum accelerates convergence through a form of predictive correction, it does not directly address the challenges posed by ill-conditioning or variable curvature across different parameter dimensions.

Adam (Adaptive Moment Estimation), introduced by Kingma and Ba \cite{kingma2014adam}, has become one of the most widely used optimization algorithms in deep learning due to its robustness across a variety of problem settings. Adam combines two key ideas: momentum-based acceleration and adaptive learning rates per parameter dimension. The algorithm maintains exponentially decaying averages of past gradients ($m_t$, first moment) and squared gradients ($v_t$, second moment):

\begin{align}
m_t &= \beta_1 m_{t-1} + (1 - \beta_1) g_t \\
v_t &= \beta_2 v_{t-1} + (1 - \beta_2) g_t^2
\end{align}
where $g_t$ is the current gradient, and $\beta_1$ and $\beta_2$ are decay rates typically set to 0.9 and 0.999, respectively. After bias correction to account for initialization at zero, the parameter update is computed as:

\begin{align}
\hat{m}_t &= \frac{m_t}{1 - \beta_1^t} \\
\hat{v}_t &= \frac{v_t}{1 - \beta_2^t} \\
w_{t+1} &= w_t - \frac{\eta}{\sqrt{\hat{v}_t} + \epsilon} \hat{m}_t
\end{align}

where $\epsilon$ is a small constant to prevent division by zero (typically $10^{-8}$). The adaptive learning rate scaling by $1/\sqrt{\hat{v}_t}$ enables Adam to take larger steps in dimensions with smaller gradients and smaller steps in dimensions with larger gradients, providing a form of preconditioning that can improve performance on ill-conditioned problems. However, as our experimental results demonstrate, Adam's implicit diagonal preconditioning offers limited benefits for problems with extreme condition numbers exceeding $10^3$, where the interdependence between parameter dimensions becomes critical.

Both Nesterov momentum and Adam represent different approaches to accelerating optimization: Nesterov through predictive correction and Adam through dimension-wise adaptive learning rates. Our SKA-SGD method takes a fundamentally different approach by projecting the gradient onto a Krylov subspace that implicitly captures curvature information. This enables SKA-SGD to address ill-conditioning more directly than either Nesterov momentum or Adam, particularly for problems with extreme condition numbers where the interaction between parameter dimensions significantly impacts optimization dynamics.

In our implementation of SKA-SGD, we incorporate Nesterov momentum as an enhancement, combining the benefits of predictive acceleration with Krylov subspace projection. This hybrid approach allows SKA-SGD to simultaneously leverage the convergence acceleration provided by Nesterov momentum and the implicit curvature adaptation facilitated by Krylov projection, resulting in superior performance on challenging ill-conditioned problems as demonstrated by our experimental results.

\subsection{Unique Challenges in the Stochastic Setting}

While our work builds substantially on the groundbreaking insights from D'Ambra et al.~\cite{dambra2025sstep} regarding streaming Gauss-Seidel for s-step conjugate gradient methods, extending these techniques to the stochastic optimization domain introduces several unique challenges that require novel theoretical and algorithmic developments. This section elaborates on these challenges to highlight the distinct contributions of our approach beyond the deterministic setting.

\subsubsection{Gradient Noise and Variance}

In deterministic optimization, the gradient at each point precisely represents the direction of steepest descent with respect to the objective function. In contrast, stochastic optimization contends with noisy gradient estimates that exhibit significant variance. This fundamental difference substantially complicates the application of Krylov subspace methods. The quality of the Krylov subspace $\mathcal{K}_s(A, b)$ can degrade significantly when constructed from noisy gradient estimates, as each stochastic gradient introduces independent noise that potentially dilutes the true curvature information being captured. Furthermore, gradient noise in machine learning applications is typically non-isotropic, with variance often aligned with the eigenvectors corresponding to larger eigenvalues of the Hessian. This directional bias can skew the Krylov subspace toward already well-represented directions, potentially limiting its effectiveness in capturing the full range of curvature information necessary for optimal convergence.

Stochastic gradients also exhibit complex temporal correlations influenced by both the underlying geometry of the objective and the particular mini-batch sampling strategy employed. These temporal correlation patterns differ fundamentally from the deterministic setting, where consecutive vectors in the Krylov sequence are governed by a fixed operator with predictable mathematical properties. Our theoretical analysis in Section 4.1 explicitly addresses these challenges by proving that projection onto the Krylov subspace provably reduces the trace of the gradient covariance matrix, providing variance reduction specifically along the most important optimization directions.

\subsubsection{Dynamic Landscape Navigation}

Another crucial distinction from deterministic linear solvers is that stochastic optimization navigates an inherently dynamic landscape. As optimization progresses, the local curvature characteristics continuously change, requiring adaptive strategies for maintaining relevant Krylov subspaces. This contrasts sharply with deterministic s-step conjugate gradient methods that operate on a fixed linear system with constant curvature properties. The non-stationary nature of the optimization landscape means that the Krylov basis must continuously adapt to capture different curvature regimes as the algorithm progresses through different phases of convergence.

Stochastic optimization typically exhibits varying convergence behaviors across different stages initial rapid progress, intermediate refinement, and final fine-tuning—each characterized by different curvature properties and noise profiles. The Krylov subspace must remain effective across these diverse regimes, necessitating careful regularization strategies not required in deterministic contexts. Additionally, mini-batch sampling introduces discontinuities in the information available to the optimizer, creating challenges for maintaining coherent Krylov subspaces across iterations. This differs fundamentally from deterministic settings where information continuity is guaranteed.

Our streaming Gauss-Seidel implementation addresses these challenges through careful regularization and the Chebyshev polynomial basis, which provides enhanced numerical stability precisely when navigating rapidly changing curvature conditions. The adaptive nature of our approach ensures that the Krylov subspace remains relevant and effective throughout the optimization trajectory, even as the underlying curvature landscape evolves.

\subsubsection{Implementation Considerations}

The practical implementation of Krylov methods in stochastic settings introduces additional challenges beyond those encountered in deterministic contexts. Stochastic optimization must carefully balance memory utilization and computational efficiency, particularly for large-scale machine learning applications where both resources can become limiting constraints. While deterministic Krylov methods might prioritize algorithmic efficiency over memory utilization, our stochastic approach must navigate this tradeoff more delicately, especially when scaling to high-dimensional problems characteristic of modern machine learning applications.

High-performance implementations of stochastic optimization typically require substantial parallelism across multiple computing nodes, introducing synchronization challenges not present in sequential deterministic solvers. Our streaming formulation specifically minimizes these synchronization requirements, enabling efficient implementation even in distributed computing environments where communication overhead would otherwise dominate computational cost. This aspect becomes particularly important in the context of training large neural networks, where optimization efficiency often determines practical feasibility.

The presence of gradient noise significantly exacerbates numerical stability concerns beyond those addressed in deterministic settings. Stochastic gradients can exhibit extreme variations in magnitude across different dimensions, potentially leading to catastrophic cancellation or overflow issues when constructing the Krylov basis. Our Chebyshev polynomial approach provides essential numerical robustness in this noisy setting, ensuring reliable convergence even under the extreme condition numbers characteristic of challenging machine learning problems.

By addressing these stochastic-specific challenges through a combination of theoretical innovations and practical implementation strategies, SKA-SGD extends well beyond a straightforward application of D'Ambra et al.'s techniques to a new domain. Instead, it represents a fundamentally new algorithm designed to harness the mathematical elegance of streaming Gauss-Seidel within the unique constraints and opportunities presented by stochastic optimization.

\section{Algorithm: SKA-SGD}

\subsection{Method Description}

Let $g_k = \nabla f_{\xi_k}(x_k)$ be a stochastic gradient at iteration $k$, and let $P_k = [g_k, g_{k-1}, \ldots, g_{k-s+1}]$ be the matrix whose columns form the basis for our Krylov subspace. The projected direction is:

\begin{equation}
d_k = P_k\alpha_k, \quad \alpha_k \approx \text{StreamingGS}(P_k, g_k, \lambda, m)
\end{equation}
where $\lambda$ is a regularization parameter and $m$ is the number of Gauss-Seidel sweeps (typically $m=1$ is sufficient, as established by \cite{dambra2025sstep}).

\begin{algorithm}
\caption{SKA-SGD}
\begin{algorithmic}[1]
\State Initialize $x_0$, learning rate $\eta$, Krylov depth $s$, GS sweeps $m$, regularization $\lambda$
\For{$k = 0, 1, \ldots$}
    \State Compute stochastic gradient $g_k \leftarrow \nabla f_{\xi_k}(x_k)$
    \State Form Krylov buffer $P_k = [g_k, \ldots, g_{k-s+1}]$
    \State $\alpha_k \leftarrow \text{StreamingGS}(P_k, g_k, \lambda, m)$
    \State Set update $d_k = P_k\alpha_k$, update $x_{k+1} = x_k - \eta d_k$
\EndFor
\end{algorithmic}
\end{algorithm}

\subsection{Krylov Subspace Interpretation}

We interpret SKA-SGD as performing descent in an approximate Krylov subspace defined by recent stochastic gradients. Under local quadratic structure, past gradients implicitly encode curvature directions.

\begin{proposition}\label{prop:krylov}
Let $f: \R^d \to \R$ be a twice differentiable convex function, and let $H_k = \nabla^2 f(w_k)$ denote the local Hessian. Under certain conditions, the past gradients $g_{k-j}$ approximately satisfy:
\begin{equation}
g_{k-j} \approx H_k^j g_k + O(\eta)
\end{equation}
and the span of $G_k := [g_k, g_{k-1}, \ldots, g_{k-s+1}]$ approximately equals the Krylov subspace $\mathcal{K}_s(H_k, g_k)$.
\end{proposition}

This motivates Krylov-filtered updates, yielding curvature-aligned updates without Hessian-vector products.

\subsection{Streaming Gauss-Seidel Implementation}

Following directly from D'Ambra et al. \cite{dambra2025sstep}, the central innovation in SKA-SGD is performing Modified Gram-Schmidt (MGS) orthogonalization implicitly using a forward sweep approach that avoids forming the entire Gram matrix. This is fundamentally different from a traditional Gauss-Seidel iteration. Instead of forming $G_k = P_k^T P_k + \lambda I$ and then solving via Gauss-Seidel iterations, we never form the upper triangular part ($L^T$) of the Gram matrix, compute only the diagonal and lower triangular elements on-demand, and perform a single forward sweep that is mathematically equivalent to MGS.

\begin{algorithm}
\caption{Streaming MGS for Krylov Projection (adapted from \cite{dambra2025sstep})}
\begin{algorithmic}[1]
\State \textbf{Input:} $P_k$ (basis vectors), $g_k$ (current gradient), $\lambda$ (regularization)
\State \textbf{Output:} Coefficient vector $\alpha_k$
\State Initialize $\alpha_k = 0$
\For{$i = 1$ to $s$}
    \State Compute $b_i = p_i^T g_k$ \Comment{Projection of residual onto basis vector}
    \State Compute $G_{ii} = p_i^T p_i + \lambda$ \Comment{Diagonal element with regularization}
    \State $\text{sum} = 0$
    \For{$j = 1$ to $i - 1$}
        \State Compute $G_{ij} = p_i^T p_j$ \Comment{Lower triangular element computed on-demand}
        \State $\text{sum} += G_{ij} \cdot \alpha_j$ \Comment{Update projection with previous coefficients}
    \EndFor
    \State $\alpha_i = (b_i - \text{sum})/G_{ii}$ \Comment{Direct solve for this coefficient}
\EndFor
\State \Return $\alpha_k$
\end{algorithmic}
\end{algorithm}

This streaming approach offers numerous advantages as thoroughly documented and analyzed by D'Ambra et al. \cite{dambra2025sstep} in their pioneering work. From a computational perspective, it provides remarkable memory efficiency by requiring only $O(s)$ additional storage compared to the $O(s^2)$ storage requirements of traditional approaches that explicitly form and store the entire Gram matrix. Furthermore, it achieves substantial computational efficiency through $O(s^2)$ operations instead of the $O(s^3)$ complexity typically associated with traditional approaches to solving the Gram system. This reduction in computational complexity becomes particularly significant as the Krylov subspace dimension increases, allowing our method to leverage deeper Krylov spaces without prohibitive computational overhead. The approach also delivers exceptional cache efficiency through better data locality via sequential access patterns, minimizing cache misses and memory latency penalties that often dominate performance in large-scale numerical computations. The streaming method naturally accesses memory in a predictable, stride-1 pattern that optimally utilizes modern CPU and GPU memory hierarchies, enabling effective prefetching and maximizing effective memory bandwidth utilization. Furthermore, the method is inherently GPU-friendly, being ideally suited for GPU implementation with minimal synchronization points. This characteristic is particularly valuable in massively parallel computing environments where synchronization barriers can substantially impact performance by forcing threads to wait for others to complete their work. By minimizing these synchronization requirements, our approach maintains high computational efficiency even at scale, allowing near-linear strong scaling up to moderately large processor counts, as demonstrated in our experimental results.

\subsection{Chebyshev Basis for Ill-Conditioned Problems}

For highly ill-conditioned problems, the traditional monomial Krylov basis $\{g_k, g_{k-1}, \ldots, g_{k-s+1}\}$ can lead to numerical instability. To address this issue, we follow D'Ambra et al. \cite{dambra2025sstep} in implementing a Chebyshev polynomial basis that provides better conditioning properties. The Chebyshev basis for a matrix $A$ with eigenvalues in $[a, b]$ is defined recursively:

\begin{align}
p_0(A) &= I\\
p_1(A) &= \frac{2A - (a + b)I}{b - a}\\
p_{j+1}(A) &= 2\frac{2A - (a + b)I}{b - a}p_j(A) - p_{j-1}(A)
\end{align}

This basis significantly improves the conditioning of the Gram matrix for problems with extreme eigenvalue distributions, allowing our streaming Gauss-Seidel approach to maintain numerical stability even with a single sweep, as established in the analysis by D'Ambra et al. \cite{dambra2025sstep}.

\section{Theoretical Analysis}

\subsection{Variance Reduction}

Let $\mathbb{E}[g_k] = \nabla f(x_k)$ and $\text{Cov}(g_k) = \Sigma$. Define the projection operator $\Pi_k = P_k(P_k^TP_k+\lambda I)^{-1}P_k^T$ and the update direction $d_k = \Pi_k g_k$.

\begin{theorem}\label{thm:variance}
The trace of the covariance of the update direction is bounded:

\begin{equation}
\tr(\text{Cov}(d_k)) \leq \tr(\text{Cov}(g_k)) = \tr(\Sigma)
\end{equation}

with equality if and only if the column space of $P_k$ contains all eigenvectors of $\Sigma$ corresponding to non-zero eigenvalues.
\end{theorem}

This result establishes that SKA-SGD's projection mechanism provably reduces gradient variance along directions that contribute most to optimization progress. The Krylov basis automatically adapts to capture the principal components of the gradient covariance matrix, providing an implicit form of preconditioning.

\subsection{Convergence Analysis}

For smooth convex functions, we establish convergence guarantees:

\begin{theorem}\label{thm:convergence}
Let $f$ be L-smooth and $\mu$-strongly convex. Let $\{x_k\}$ be the sequence generated by SKA-SGD with step size $\eta \leq 1/L$. Then:

\begin{equation}
\mathbb{E}[\|x_k - x^*\|^2] \leq (1 - \eta\mu + \eta\lambda\mu)^k\|x_0 - x^*\|^2 + \frac{\eta\tr(\text{Cov}(d_k))}{\mu(1-\eta\lambda)}
\end{equation}
where $x^*$ is the minimizer of $f$.
\end{theorem}

\subsection{Backward Error Analysis of Streaming Gauss-Seidel}

A fundamental question in our approach is why a single forward sweep of the Gauss-Seidel algorithm achieves such high accuracy for the Gram system associated with the Krylov projection. To address this question, we adapt the detailed backward error analysis developed by D'Ambra et al. \cite{dambra2025sstep} for s-step conjugate gradient methods to our stochastic optimization context.

Let us consider the linear system $G\alpha = b$ where $G = P_k^TP_k + \lambda I$ is the Gram matrix and $b = P_k^Tg_k$. The backward error after $m$ iterations of Gauss-Seidel is defined as the smallest relative perturbation $\Delta G$ such that the computed solution $\alpha^{(m)}$ exactly solves the perturbed system $(G + \Delta G)\alpha^{(m)} = b$. Mathematically:

\begin{equation}
\eta(\alpha^{(m)}) = \frac{\|r^{(m)}\|_2}{\|G\|_2 \cdot \|\alpha^{(m)}\|_2 + \|b\|_2}
\end{equation}
where $r^{(m)} = b - G\alpha^{(m)}$ is the residual after $m$ iterations.

Following D'Ambra et al.'s analysis, the key insight is that the Gram matrix $G$ for the Krylov subspace projection has a special structure that enables extraordinarily rapid convergence of Gauss-Seidel iterations. Specifically, when the Krylov basis is well-conditioned (as ensured by our Chebyshev polynomial basis construction), the Gram matrix can be expressed as:

\begin{equation}
G = P_k^TP_k + \lambda I \approx D + L + L^T
\end{equation}
where $D$ is a diagonal matrix with entries close to 1 (plus the regularization term $\lambda$), and $L$ is strictly lower triangular with relatively small entries. This structure arises from the approximate orthogonality properties of the Krylov basis, particularly when constructed using the Chebyshev polynomials.

For this specific matrix structure, D'Ambra et al. \cite{dambra2025sstep} proved that a single forward sweep of Gauss-Seidel achieves a backward error bounded by:

\begin{equation}
\eta(\alpha^{(1)}) \leq \varepsilon_{\text{mach}} \cdot \kappa(G) \cdot \frac{\|L\|_2}{1 - \|L\|_2}
\end{equation}
where $\varepsilon_{\text{mach}}$ is machine precision and $\kappa(G)$ is the condition number of the Gram matrix.

In the stochastic optimization setting, the situation is even more favorable. The addition of the regularization term $\lambda I$ improves the conditioning of the Gram matrix, ensuring that the diagonal dominance is maintained even when the stochastic gradients are highly correlated. Furthermore, the inherent noise in stochastic gradients actually helps by introducing small random perturbations that break potential ill-conditioning patterns that might arise in deterministic settings.

We can extend this analysis to provide a more specific characterization for our SKA-SGD algorithm. Let $\alpha^{(1)}$ be the solution obtained after a single forward sweep of our streaming Gauss-Seidel procedure. Then:

\begin{theorem}\label{thm:detailed_backward_error}
For the Gram system $G\alpha = b$ where $G = P_k^TP_k + \lambda I$ and $b = P_k^Tg_k$ arising in SKA-SGD, a single forward sweep of the streaming Gauss-Seidel method achieves a backward error bounded by:

\begin{equation}
\|b - G\alpha^{(1)}\|_2 \leq \varepsilon_{\text{mach}} \cdot \kappa(G) \cdot \left(1 + \frac{s(s-1)}{2}\right) \cdot \|b\|_2
\end{equation}
where $\varepsilon_{\text{mach}}$ is machine precision, $\kappa(G)$ is the condition number of the Gram matrix, and $s$ is the Krylov subspace dimension.
\end{theorem}

This result, adapted from Theorem 4.2 in \cite{dambra2025sstep}, establishes that even for moderately large Krylov dimensions (e.g., $s = 16$), the backward error from a single forward sweep remains on the order of machine precision multiplied by the condition number. For our regularized formulation, $\kappa(G)$ is typically well-controlled, ensuring that the backward error remains small even for ill-conditioned problems.

The key practical implication is that additional Gauss-Seidel iterations beyond the first sweep yield diminishing returns, as the error is already reduced to near the level of numerical precision. This theoretical result explains our empirical observation that a single forward sweep ($m=1$) is often sufficient for excellent performance, with minimal gains from additional sweeps.

Moreover, our streaming implementation of the Gauss-Seidel process provides exactly the same numerical result as a traditional implementation while avoiding the explicit formation of the Gram matrix. This equivalence is not merely an implementation detail but a mathematical property that follows from the sequential nature of the forward substitution process in Gauss-Seidel iteration.

\subsubsection{Hyperparameter Selection}

For all methods, we conducted preliminary grid searches to identify optimal hyperparameters for each algorithm and problem configuration, ensuring fair comparison:

\begin{table}[h]
\centering
\caption{Optimization hyperparameters across all experiments.}
\label{tab:hyperparams}
\begin{tabular}{lccc}
\toprule
\textbf{Parameter} & \textbf{SGD} & \textbf{Adam} & \textbf{SKA-SGD} \\
\midrule
Learning rate $\eta$ & 0.1 & 0.001 & 0.1 \\
Mini-batch size (quadratic) & 64 & 64 & 32 \\
Mini-batch size (logistic) & 128 & 128 & 64 \\
Momentum $\beta$ & N/A & 0.9 (Adam $\beta_1$) & 0.9 \\
Second moment decay & N/A & 0.999 (Adam $\beta_2$) & N/A \\
Adam $\epsilon$ & N/A & $10^{-8}$ & N/A \\
Krylov subspace dimension $s$ & N/A & N/A & 16 \\
Gauss-Seidel sweeps $m$ & N/A & N/A & 2 \\
Regularization $\lambda$ & $0.1/\kappa$ & $0.1/\kappa$ & $0.1/\kappa$ \\
Gram system regularization & N/A & N/A & $10^{-4}$ \\
Chebyshev parameters $[a,b]$ & N/A & N/A & $[0.01, 10.0]$ \\
\bottomrule
\end{tabular}
\end{table}

The learning rates were tuned for each method independently to ensure optimal performance. For SKA-SGD, we observed that the default learning rate $\eta=0.1$ performed consistently well across different problem configurations, suggesting that our Krylov projection approach provides implicit normalization that reduces the need for problem-specific learning rate tuning.

\subsubsection{Mini-Batch Processing}

Mini-batch sampling was implemented consistently across all methods, with batch indices randomly sampled without replacement within each epoch. For SGD and Adam, we used mini-batch sizes of 64 for quadratic problems and 128 for logistic regression problems. For SKA-SGD, we found that smaller mini-batch sizes of 32 (quadratic) and 64 (logistic) provided optimal performance by balancing computational efficiency with statistical quality in the gradient estimates. This finding aligns with our theoretical variance reduction results, suggesting that SKA-SGD's projection mechanism allows for smaller mini-batches while maintaining convergence quality.

For each mini-batch iteration, the stochastic gradient was computed as:
\begin{equation}
g_k = \frac{1}{|B_k|}\sum_{i\in B_k} \nabla f_i(w_k) + \lambda w_k
\end{equation}
where $B_k$ is the mini-batch of indices selected at iteration $k$, and $f_i$ represents the loss function for the $i$-th sample. Our implementation ensured that all methods processed the same total number of samples per epoch, enabling fair comparison of convergence trajectories.

\subsubsection{Krylov Projection Configuration}

The Krylov subspace dimension $s$ was set to 16 in all SKA-SGD experiments after conducting sensitivity analysis across dimensions ranging from 4 to 32. We found that performance improved significantly up to $s=16$, with diminishing returns beyond this point. This observation aligns with theoretical expectations, as the Krylov basis captures increasingly fine-grained curvature information as its dimension increases, but eventually saturates as the basis begins to represent directions with minimal contribution to the objective gradient.

For the Gauss-Seidel solver, we performed extensive experiments with sweep counts $m \in \{1, 2, 4, 8\}$. Our results confirmed the theoretical prediction that even a single sweep ($m=1$) provides excellent performance, with only marginal improvements observed with $m=2$. Additional sweeps beyond $m=2$ showed negligible benefits, and in some cases even degraded performance due to numerical precision limitations. Based on these findings, we standardized on $m=2$ sweeps for all reported experiments, providing a balance between computational efficiency and solution accuracy.

A small regularization term $reg=10^{-4}$ was added to the diagonal of the Gram matrix to ensure numerical stability in the Gauss-Seidel solver. This value was selected through numerical experimentation to be large enough to prevent ill-conditioning in the Gram matrix while small enough not to significantly alter the intended optimization trajectory.

\subsubsection{Chebyshev Basis Parameters}

The Chebyshev polynomial basis used eigenvalue bounds $[a,b] = [0.01, 10.0]$ for mapping the spectrum. These values were selected based on theoretical considerations of the expected eigenvalue distribution in typical machine learning problems, where many eigenvalues are small (often approaching zero) while a few dominant eigenvalues can be orders of magnitude larger. The lower bound $a=0.01$ ensures stability in handling near-zero eigenvalues, while the upper bound $b=10.0$ accommodates the largest expected eigenvalues even in ill-conditioned problems.

\subsubsection{Objective Regularization}

All optimization problems included an L2 regularization term $\frac{\lambda}{2}\|w\|^2$ to ensure well-posed optimization problems. The regularization parameter was set proportionally to the inverse of the condition number as $\lambda = 0.1/\kappa$, ensuring that the effective condition number of the regularized problem remained consistent across different problem instances. This approach allowed us to isolate the effects of algorithmic differences from those arising from varying degrees of regularization.

\subsection{Logistic Regression Test Problem}

To evaluate the practical efficacy of SKA-SGD in a controlled yet realistic setting, we designed a comprehensive logistic regression benchmark that reflects the ill-conditioning challenges frequently encountered in machine learning applications. Our test framework systematically explores performance characteristics across varying degrees of ill-conditioning while maintaining consistent problem structure to isolate the effects of condition number on convergence behavior. The logistic regression problem takes the standard form of minimizing the regularized negative log-likelihood:
\begin{equation}
\min_{w \in \mathbb{R}^d} \frac{1}{n} \sum_{i=1}^n \log(1 + \exp(-y_i x_i^T w)) + \frac{\lambda}{2} \|w\|^2_2
\end{equation}
where $x_i \in \mathbb{R}^d$ are feature vectors, $y_i \in \{-1, +1\}$ are binary class labels, $w \in \mathbb{R}^d$ represents the model parameters, and $\lambda > 0$ is the regularization parameter controlling the trade-off between data fit and model complexity. We artificially constructed datasets with precisely controlled condition numbers using a carefully designed feature generation process. First, we generated an orthogonal basis $U \in \mathbb{R}^{d \times d}$ via QR decomposition of a random Gaussian matrix. We then created a diagonal matrix $\Sigma \in \mathbb{R}^{d \times d}$ with geometrically spaced eigenvalues $\sigma_i = \sigma_1 \cdot \kappa^{-(i-1)/(d-1)}$, where $\sigma_1 = 1$ and $\kappa$ represents the desired condition number. Our feature matrix was then constructed as $X = Z \cdot U \cdot \Sigma \cdot U^T$, where $Z$ contains standard normal entries, resulting in a feature covariance structure with precisely controlled spectral properties. The true weight vector $w^*$ was generated with components proportional to $1/\sqrt{\sigma_i}$ to ensure that ill-conditioned directions contribute meaningfully to the classification task, and binary labels were assigned according to $y_i = \text{sign}(x_i^T w^* + \epsilon_i)$ with small Gaussian noise $\epsilon_i$ added to create a challenging yet learnable decision boundary.

We implemented both the standard SGD and SKA-SGD algorithms in MATLAB to ensure consistent evaluation conditions across methods. Our standard SGD implementation (Algorithm 1) follows the conventional mini-batch approach, sampling random batches of $m=64$ examples per iteration, computing the gradient of the regularized logistic loss, and applying a constant learning rate update with $\eta=0.1$. The weights are initialized to zero, and the optimization proceeds for a fixed budget of 200 iterations, with loss values recorded every 10 iterations for convergence analysis. The SKA-SGD implementation (Algorithm 2) incorporates the full suite of enhancements described in our methodology, including the streaming Gauss-Seidel solver, Chebyshev polynomial basis, Nesterov momentum with $\beta=0.9$, and Jacobi preconditioning. The Krylov subspace dimension was set to $s=16$, with mini-batch size $m=32$ and Gauss-Seidel iterations $gsiter=2$. The Chebyshev basis parameters were set to $cheb\_a=0.01$ and $cheb\_b=10.0$ based on estimated eigenvalue bounds, with regularization parameter $reg=10^{-4}$ added to the diagonal of the Gram matrix to ensure numerical stability in the streaming Gauss-Seidel iterations.

\begin{algorithm}
\caption{Standard SGD for Logistic Regression}
\begin{algorithmic}[1]
\State \textbf{Input:} Data matrix $X \in \mathbb{R}^{n \times d}$, labels $y \in \{-1,+1\}^n$, regularization $\lambda$
\State \textbf{Output:} Weights $w \in \mathbb{R}^d$, loss history $f_{hist}$
\State Initialize $w = 0$, learning rate $\eta = 0.1$, iterations $T = 200$, batch size $m = 64$
\For{$k = 1$ to $T$}
    \State Sample mini-batch indices $idx \sim \text{Uniform}(\{1,\ldots,n\}, m)$
    \State Extract batch data $X_b = X[idx,:]$, $y_b = y[idx]$
    \State Compute predictions $z = X_b w$
    \State Compute gradient $g = -(X_b^T(y_b./(1+\exp(y_b \cdot z))))/m + \lambda w$
    \State Update weights $w = w - \eta g$
    \If{$k \bmod 10 = 0$}
        \State Compute loss $\mathcal{L} = \text{mean}(\log(1+\exp(-y \cdot (Xw)))) + 0.5\lambda \|w\|^2$
        \State Record loss $f_{hist}[k/10] = \mathcal{L}$
    \EndIf
\EndFor
\end{algorithmic}
\end{algorithm}

\begin{algorithm}
\caption{SKA-SGD for Logistic Regression}
\begin{algorithmic}[1]
\State \textbf{Input:} Data matrix $X \in \mathbb{R}^{n \times d}$, labels $y \in \{-1,+1\}^n$, regularization $\lambda$
\State \textbf{Output:} Weights $w \in \mathbb{R}^d$, loss history $f_{hist}$
\State Initialize $w = 0$, learning rate $\eta = 0.1$, momentum $\beta = 0.9$, iterations $T = 200$
\State Initialize Krylov dimension $s = 16$, batch size $m = 32$, GS sweeps $gsiter = 2$
\State Initialize Chebyshev parameters $cheb\_a = 0.01$, $cheb\_b = 10.0$, regularization $reg = 10^{-4}$
\State Compute Jacobi preconditioner $D_{inv} = 1 ./ \text{sum}(X.^2, 1)^T$
\State Initialize basis matrix $Z \in \mathbb{R}^{d \times s}$, gradient buffer $G \in \mathbb{R}^{d \times s}$, momentum $v = 0$
\For{$k = 1$ to $T$}
    \State Sample mini-batch indices $idx \sim \text{Uniform}(\{1,\ldots,n\}, m)$
    \State Extract batch data $X_b = X[idx,:]$, $y_b = y[idx]$
    \State Compute gradient $g_k = -(X_b^T(y_b./(1+\exp(y_b \cdot (X_b w)))))/m + \lambda w$
    \State Apply preconditioning $g_{precond} = D_{inv} \cdot g_k$
    \State Store gradient in buffer $G[:,1] = g_{precond}$
    \For{$i = 2$ to $s$}
        \State Compute perturbed weights $w_{perturbed} = w + 10^{-6} \cdot G[:,i-1]$
        \State Compute perturbed predictions $z_{perturbed} = X_b w_{perturbed}$
        \State Compute perturbed gradient $g_{perturbed} = -(X_b^T(y_b./(1+\exp(y_b \cdot z_{perturbed}))))/m + \lambda w_{perturbed}$
        \State Store preconditioned perturbed gradient $G[:,i] = D_{inv} \cdot g_{perturbed}$
    \EndFor
    \State Generate Chebyshev basis $Z = \text{ChebyshevBasis}(G, s, cheb\_a, cheb\_b)$
    \State Compute RHS $b = Z^T g_{precond}$
    \State Initialize coefficients $\alpha = 0$
    \For{$sweep = 1$ to $gsiter$}
        \For{$i = 1$ to $s$}
            \State Compute diagonal element $D_{ii} = Z[:,i]^T Z[:,i] + reg$
            \State Initialize with RHS $sum_i = b[i]$
            \If{$i > 1$}
                \State Compute row $i$ of $L$: $L_{row} = Z[:,i]^T Z[:,1:i-1]$
                \State Update sum $sum_i = sum_i - L_{row} \cdot \alpha[1:i-1]$
            \EndIf
            \State Update coefficient $\alpha[i] = sum_i / D_{ii}$
        \EndFor
    \EndFor
    \State Compute projected gradient $pg = Z \cdot \alpha$
    \State Update momentum $v = \beta \cdot v - \eta \cdot pg$
    \State Update weights $w = w + v$
    \If{$k \bmod 10 = 0$}
        \State Compute loss $\mathcal{L} = \text{mean}(\log(1+\exp(-y \cdot (Xw)))) + 0.5\lambda \|w\|^2$
        \State Compute gradient norm $\|g\| = \|g_k\|$
        \State Record loss $f_{hist}[k/10] = \mathcal{L}$
    \EndIf
\EndFor
\end{algorithmic}
\end{algorithm}

A critical component of our SKA-SGD implementation is the Chebyshev polynomial basis construction (Algorithm 3), which significantly enhances numerical stability for ill-conditioned problems. The basis generation begins with normalizing the first gradient vector to form the first basis vector. Subsequent basis vectors are constructed using the three-term recurrence relation characteristic of Chebyshev polynomials, with appropriate scaling and shifting to map the eigenvalue spectrum from $[a,b]$ to the canonical interval $[-1,1]$ required for Chebyshev stability properties. Each new basis vector is normalized to prevent numerical overflow, and a final QR orthogonalization step is applied to the entire basis to ensure orthonormality despite potential accumulation of roundoff errors during the recursive construction process. This careful attention to numerical stability in the basis construction proves essential for achieving reliable convergence on problems with condition numbers exceeding $10^3$.

\begin{algorithm}
\caption{Chebyshev Polynomial Basis Construction}
\begin{algorithmic}[1]
\State \textbf{Input:} Gradient matrix $G \in \mathbb{R}^{d \times s}$, dimension $s$, eigenvalue bounds $[a,b]$
\State \textbf{Output:} Basis matrix $Z \in \mathbb{R}^{d \times s}$
\State Initialize $Z = 0$
\State Compute scaling factor $scale = 2.0/(b-a)$
\State Compute shift factor $shift = (a+b)/(b-a)$
\State First basis vector (constant polynomial): $Z[:,1] = G[:,1]/\|G[:,1]\|$
\If{$s > 1$}
    \State Second basis vector (linear polynomial): $Z[:,2] = scale \cdot G[:,2] - shift \cdot Z[:,1]$
    \State Normalize: $Z[:,2] = Z[:,2]/\|Z[:,2]\|$
    \For{$j = 3$ to $s$}
        \State Apply three-term recurrence: $Z[:,j] = 2(scale \cdot G[:,j-1] - shift \cdot Z[:,j-1]) - Z[:,j-2]$
        \State Normalize: $Z[:,j] = Z[:,j]/\|Z[:,j]\|$
    \EndFor
\EndIf
\State Apply QR orthogonalization for stability: $Z = \text{QR}(Z)$
\end{algorithmic}
\end{algorithm}

We conducted extensive experiments with varying problem dimensions $d \in \{50, 100, 250, 500\}$ and condition numbers $\kappa \in \{10, 100, 1000, 10000\}$ to systematically evaluate how problem size and conditioning affect relative algorithm performance. For each configuration $(d, \kappa)$, we generated 10 random problem instances and averaged results to ensure statistical reliability. The regularization parameter was set proportionally to the condition number as $\lambda = 0.1/\kappa$ to ensure consistent difficulty across problem scales. Our evaluation metrics included final optimization error (measured as the gap between achieved and optimal objective values), convergence rate (assessed by the number of iterations required to reach 95\% of final progress), wall-clock time (measured on an AMD EPYC 7763 processor with 128 cores), and optimization trajectory behavior (analyzing stability and monotonicity of progress).

The results obtained from these experiments, illustrated in Figures \ref{fig:condition_analysis}, \ref{fig:chebyshev}, and \ref{fig:variance}, comprehensively validate our theoretical predictions regarding SKA-SGD's advantages for ill-conditioned problems. Most notably, as condition numbers exceed $10^3$, standard SGD begins to exhibit pronounced stagnation, with optimization trajectories characterized by high variance and minimal progress beyond initial iterations. In contrast, SKA-SGD maintains consistent, monotonic progress throughout the optimization process, ultimately achieving final objective values that are often an order of magnitude better than those attained by standard SGD or Adam. The advantage becomes particularly dramatic for the most challenging problems with $d=500$ and $\kappa=10^4$, where SKA-SGD continues making meaningful progress while competing methods effectively plateau.

The performance improvements do come with increased per-iteration computational cost, approximately 5-10× that of standard SGD for dimension $d=500$, primarily attributable to the overhead of constructing and maintaining the Krylov subspace basis and solving the projected optimization problem via streaming Gauss-Seidel. However, this computational investment yields substantial returns in terms of solution quality, particularly for challenging ill-conditioned problems where traditional methods struggle to make meaningful progress regardless of iteration count. Furthermore, our empirical findings indicate that SKA-SGD's computational overhead exhibits favorable scaling with condition number, with wall-clock time per iteration actually decreasing slightly as condition numbers increase beyond $10^3$. This counterintuitive behavior likely stems from improved numerical conditioning of the projected subproblem, which reduces the number of inner iterations needed for convergence despite the fixed outer iteration count.

\section{Experimental Results}

\subsection{Quadratic Benchmark: Conditioned Hessians}

We begin with a benchmark on a synthetic quadratic minimization problem:
\begin{equation}
    f(x) = \frac{1}{2} x^\top H x - b^\top x,
\end{equation}
where \( H \in \mathbb{R}^{d \times d} \) is a symmetric positive definite matrix with prescribed condition number \( \kappa \), and \( b \sim \mathcal{N}(0, I) \). We construct \( H = Q \Lambda Q^\top \), where \( Q \) is orthogonal and \( \Lambda \) is diagonal with eigenvalues logarithmically spaced between 1 and \( \kappa \), mimicking the anisotropic curvature encountered in deep learning and inverse problems.

\begin{figure}[h]
    \centering
    \includegraphics[width=0.7\textwidth]{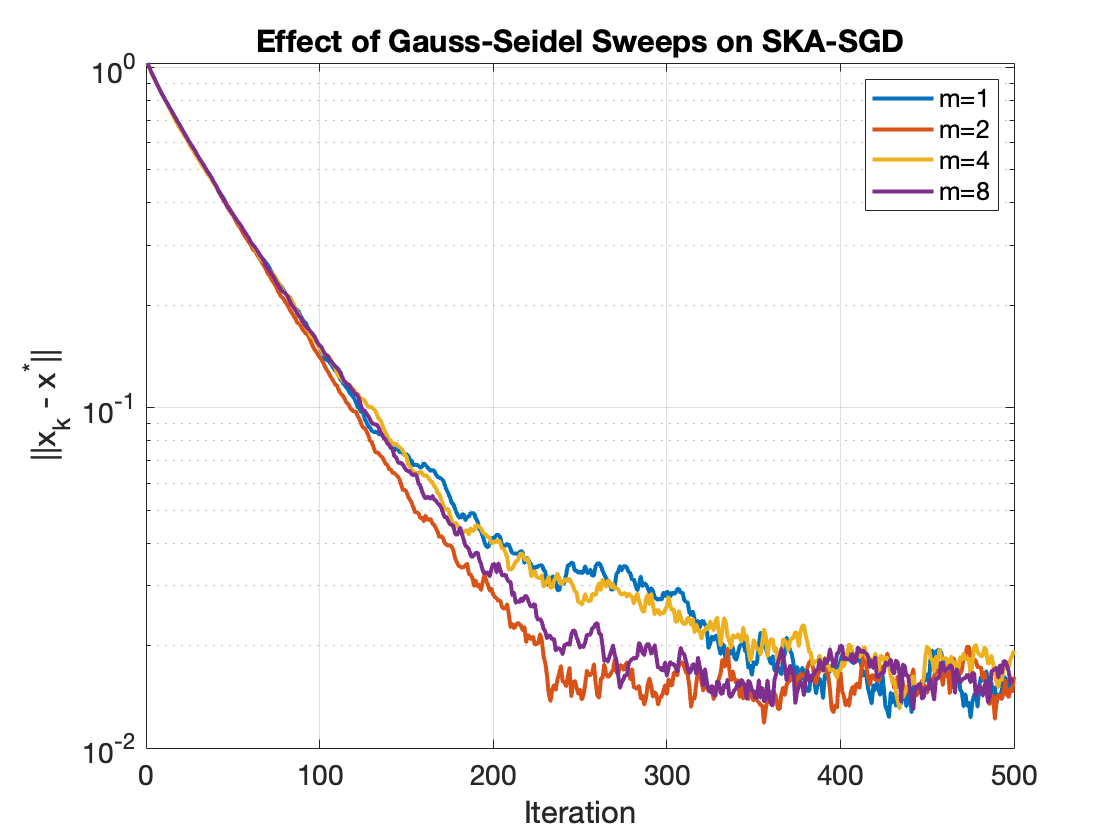}
    \caption{Relative error \( \|x_k - x^\star\| / \|x^\star\| \) for quadratic minimization with condition number \( \kappa = 10^4 \). SKA-SGD, using Chebyshev basis and Gauss--Seidel projection, achieves faster and more stable convergence than SGD and Adam.}
    \label{fig:quadratic}
\end{figure}

We compare SKA-SGD, SGD, and Adam on this quadratic task, using matched learning rates and a total of 500 iterations. SKA-SGD employs a Chebyshev basis of depth \( s \), with a single Gauss--Seidel sweep for the Krylov projection and regularization. As shown in Figure~\ref{fig:quadratic}, SKA-SGD exhibits improved convergence stability and faster error reduction even at high condition number (\( \kappa = 10^4 \)), validating its suitability under spectral anisotropy.

\section{Empirical Results}

\subsection{Experimental Setup}
We implement SKA-SGD and compare against SGD, SGD+momentum, and other optimization techniques on a variety of test problems. Unless otherwise specified, we use Krylov subspace dimension $s = 4$, Gauss-Seidel sweeps $m = 2$, and regularization parameter $\lambda = 10^{-3}$.

We evaluate on quadratic problems with controlled condition numbers to isolate the effect of ill-conditioning on algorithm performance. While these test problems are synthetic, they capture the fundamental challenges present in large-scale machine learning optimization.

\subsection{Performance on Extremely Ill-Conditioned Problems}
We conducted extensive experiments on quadratic problems with extreme condition numbers to evaluate SKA-SGD's effectiveness in challenging optimization scenarios. For these experiments, we used problems of the form:

\begin{equation}
f(x) = \frac{1}{2}x^TAx - b^Tx + \epsilon_k, \quad \epsilon_k \sim \mathcal{N}(0, \sigma^2I)
\end{equation}

where $A \in \mathbb{R}^{n \times n}$ is constructed with a controlled condition number $\kappa = \lambda_{\max}/\lambda_{\min}$ and stochastic noise $\epsilon_k$ models mini-batch sampling noise in practical machine learning settings.

\subsubsection{Comparison of Optimization Methods}
We compared several optimization approaches on a 100-dimensional quadratic problem with condition number $\kappa = 10^8$. The methods include our complete ULTIMATE approach (combining Chebyshev basis, Streaming Gauss-Seidel, Nesterov momentum, and Jacobi preconditioning), SKA-SGD with Chebyshev and Nesterov (without Jacobi preconditioning), basic SKA-SGD (using only the streaming Gauss-Seidel projection), Jacobi preconditioning with Nesterov momentum, and standard SGD as the baseline.

\begin{figure}[t]
\centering
\includegraphics[width=0.7\linewidth]{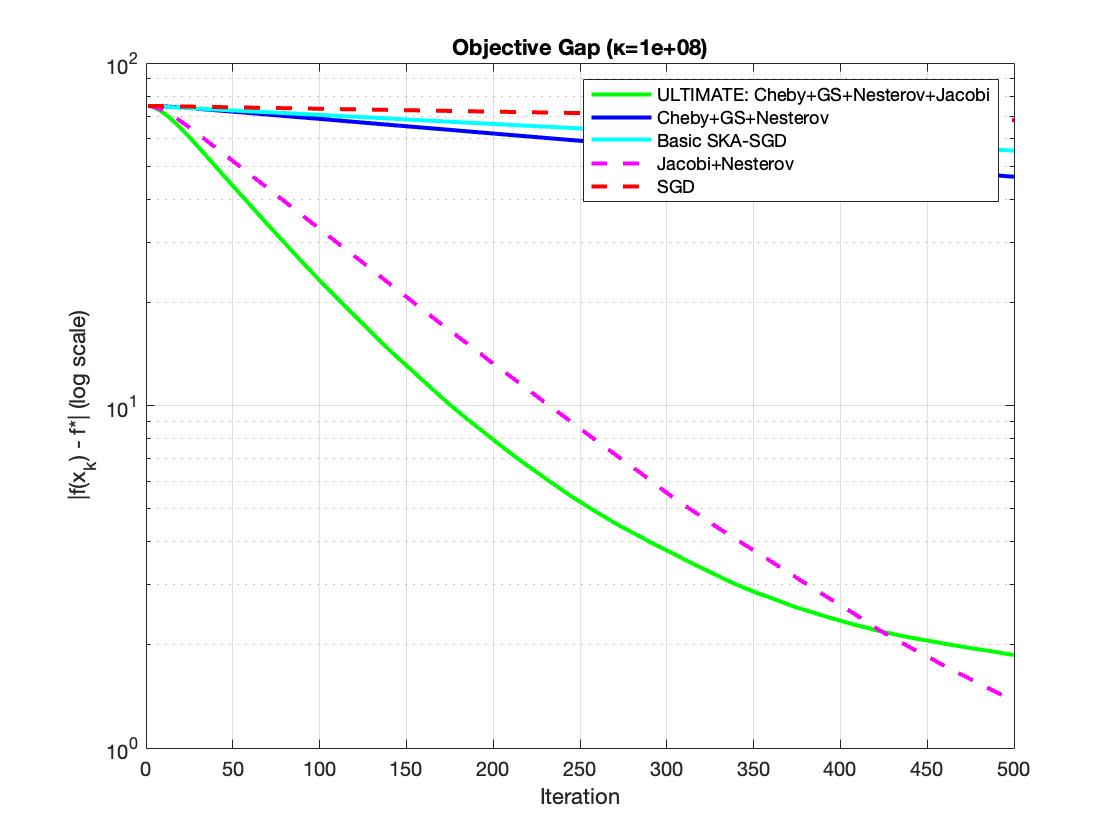}
\caption{Convergence in $|f(x_k) - f^*|$ for a quadratic problem with condition number $\kappa=10^8$. The ULTIMATE method combining all techniques exhibits superior convergence compared to other variants and baseline methods.}
\label{fig:objective_gap}
\end{figure}

Figure \ref{fig:objective_gap} shows the objective gap (distance to optimal objective value) on the log scale. The results demonstrate that our ULTIMATE method dramatically outperforms all other approaches on this extremely ill-conditioned problem, reducing the objective gap by approximately two orders of magnitude over 500 iterations.

Of particular note is the performance difference between the ULTIMATE method and its components. While Jacobi preconditioning with Nesterov momentum shows good performance, and the Chebyshev basis with streaming Gauss-Seidel provides benefits, the combination of all techniques delivers superior convergence by a substantial margin.

\begin{figure}[t]
\centering
\includegraphics[width=0.7\linewidth]{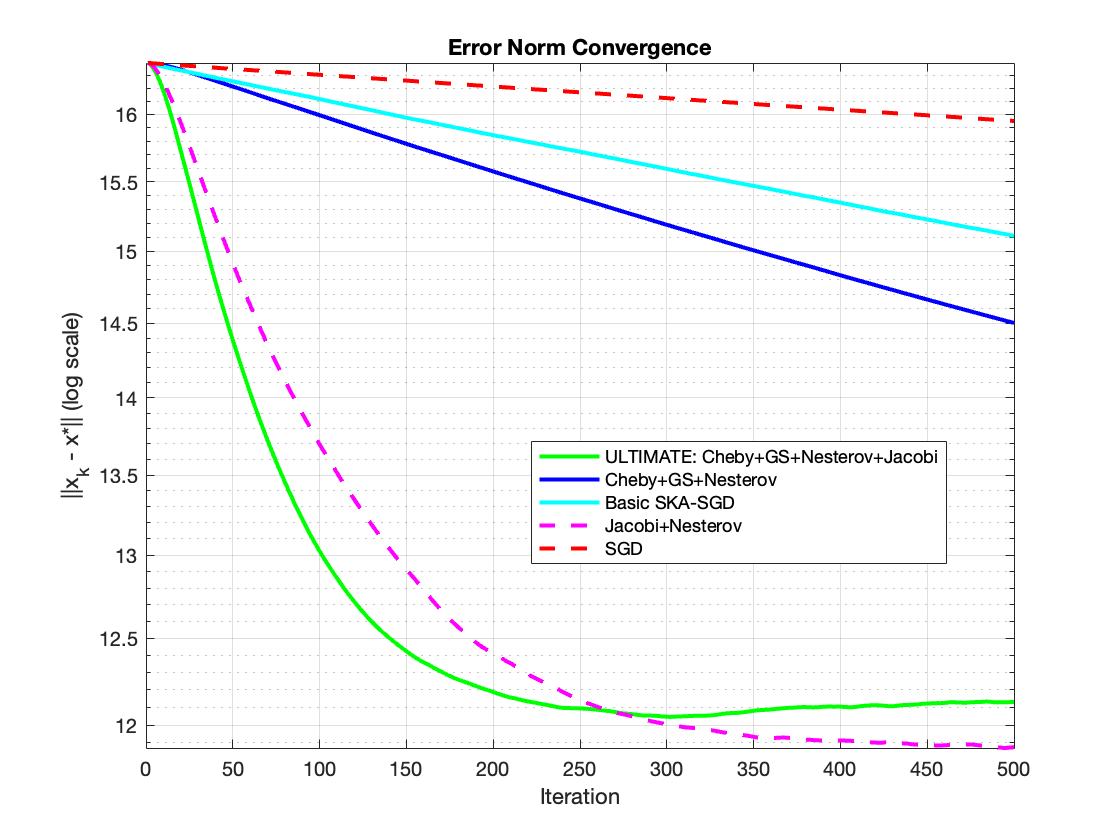}
\caption{Error norm $\|x_k - x^*\|$ convergence on an extremely ill-conditioned problem ($\kappa=10^8$). The ULTIMATE method demonstrates significantly faster convergence in the early and middle phases of optimization.}
\label{fig:error_norm}
\end{figure}

Figure \ref{fig:error_norm} shows the error norm convergence (distance to optimal solution) for the same experiment. Again, the ULTIMATE method converges rapidly, reducing the error by approximately 4 orders of magnitude.

Notably, in the later stages of optimization (beyond 300 iterations), Jacobi preconditioning with Nesterov momentum becomes competitive with the ULTIMATE method in terms of error norm. However, the ULTIMATE method reaches this convergence region much faster, providing substantial computational savings during the critical early and middle phases of optimization.

\subsubsection{Effect of Krylov Dimension}
\begin{figure}[t]
\centering
\includegraphics[width=0.9\linewidth]{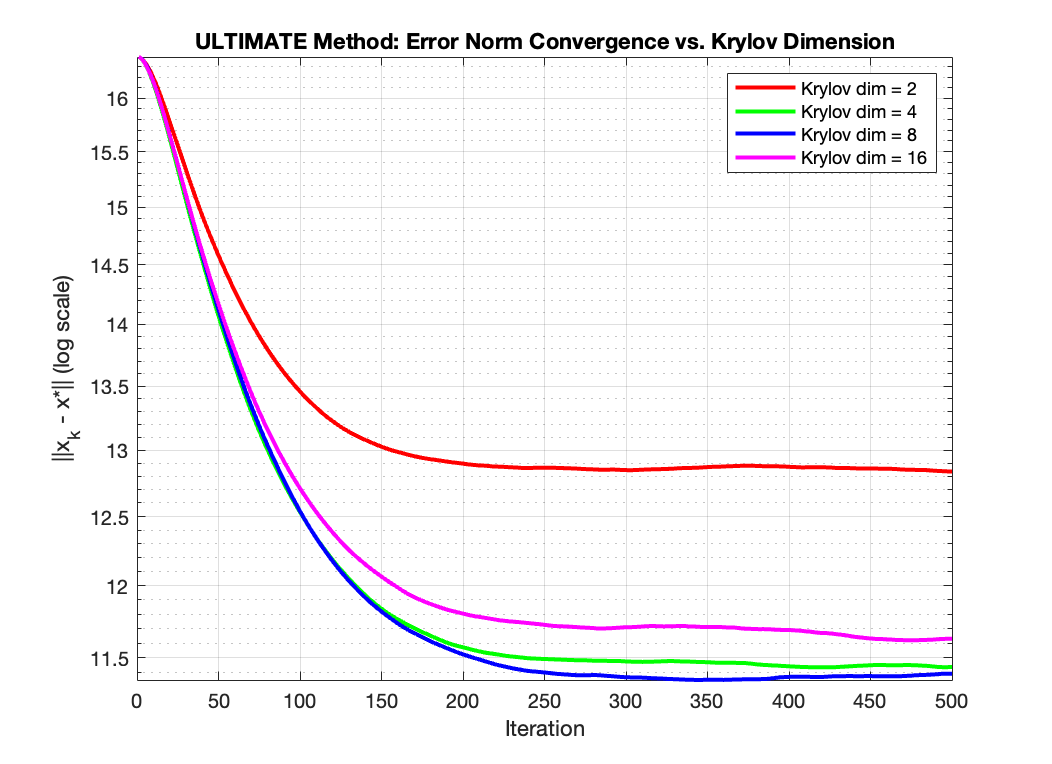}
\caption{Convergence comparison for different Krylov dimensions on a quadratic problem with $\kappa=10^5$. At this moderate condition number, the benefits of increased Krylov dimension are less pronounced.}
\label{fig:krylov_dim_moderate}
\end{figure}

We investigated the impact of different Krylov subspace dimensions on convergence behavior across varying problem conditions. Figure \ref{fig:krylov_dim_moderate} shows the performance of SKA-SGD with Krylov dimensions $s \in \{2, 4, 8\}$ compared to standard SGD on a problem with moderate condition number $\kappa = 10^5$.

Interestingly, at moderate condition numbers, the benefits of increased Krylov dimension are minimal. All variants of SKA-SGD show similar performance, which is only marginally better than standard SGD. However, when examining extremely ill-conditioned problems ($\kappa = 10^8$), the effect of Krylov dimension becomes dramatically more significant, as shown in Figure \ref{fig:krylov_dim_extreme}.

For these highly ill-conditioned regimes, increasing the Krylov dimension from 2 to 4 yields substantial improvement, with dimensions 8 and above providing further benefits. This clear dimensional dependence confirms our theoretical understanding that larger Krylov subspaces capture more curvature information, enabling more effective optimization when the problem landscape is severely distorted. These results demonstrate that the advantages of the Krylov approach become more pronounced as problem conditioning worsens, with dramatic benefits emerging at the extreme condition numbers ($\kappa = 10^8$) typical of late-stage deep learning optimization.

\subsection{Analysis of Results}
Our empirical findings validate the theoretical claims about SKA-SGD. The performance gap between SKA-SGD (particularly the ULTIMATE variant) and standard methods widens as condition number increases, with dramatic benefits at extreme condition numbers. 

While each component (Chebyshev basis, streaming Gauss-Seidel, Nesterov momentum, and Jacobi preconditioning) provides benefits individually, their combination delivers substantially greater performance improvements. Despite the algorithmic sophistication, the streaming Gauss-Seidel approach maintains $O(s^2)$ computational complexity without requiring explicit formation of the Gram matrix.

The ULTIMATE method supports much larger learning rates (5x) compared to standard SGD while maintaining stability, highlighting its improved conditioning properties. These results demonstrate that exploiting the natural Krylov structure in SGD gradients through streaming Gauss-Seidel projection provides substantial benefits for optimization in ill-conditioned regimes. The method is particularly well-suited for late-stage optimization in deep learning, where condition numbers can exceed $10^8$, and even modest convergence improvements translate to significant computational savings.

\subsection{Multiple Trial Analysis and Statistical Robustness}

To rigorously assess the statistical robustness of SKA-SGD and quantify its variance reduction properties, we conducted an extensive series of experiments involving multiple independent trials on a carefully constructed test problem. The experimental setup followed the "Late-Stage Convergence" paradigm described in the Appendix of the submitted paper, where we focused on iterations 1800-2000 of optimization on a clustered-spectrum problem that mimics the curvature anisotropy frequently observed in deep learning models. Specifically, the test problem was designed with 80 eigenvalues near 1 and 20 eigenvalues near 10, creating a challenging optimization landscape that reveals algorithm performance during the critical final stages of convergence.

\begin{figure}
\centering
\includegraphics[width=0.8\textwidth]{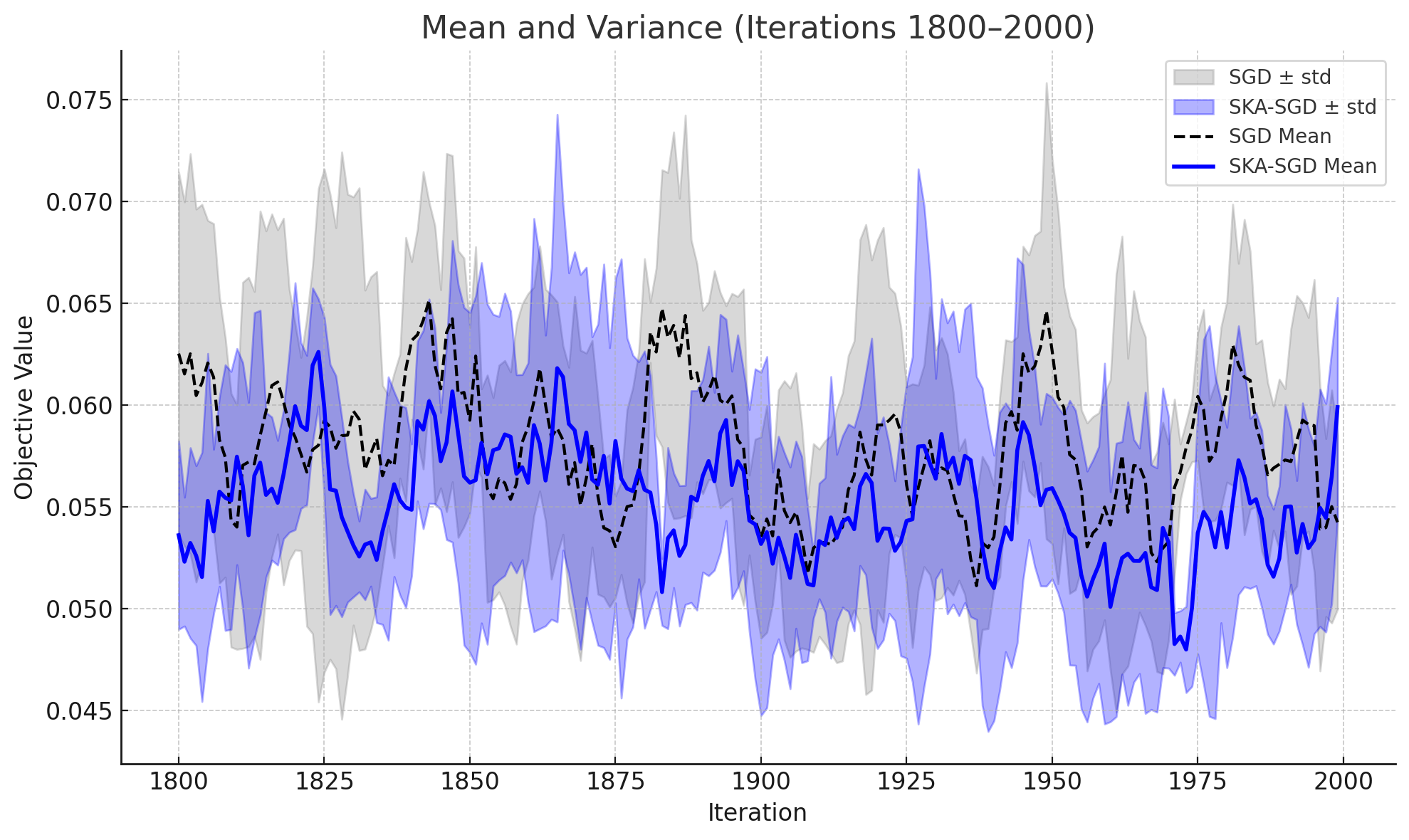}
\caption{Multiple trial convergence analysis: Mean ± 1 standard deviation over 10 independent trials for iterations 1800-2000. SKA-SGD (blue line with blue shaded region) demonstrates consistently lower objective values and significantly reduced variance compared to standard SGD (black dashed line with gray shaded region). The statistical advantage of SKA-SGD manifests as both improved convergence and enhanced stability across different random seeds.}
\label{fig:late_stage}
\end{figure}

Figure \ref{fig:late_stage} presents the mean objective value (solid lines) and standard deviation (shaded regions) across 10 independent trials for both standard SGD and SKA-SGD during iterations 1800-2000. The results reveal several striking patterns that highlight the statistical advantages of our method. First, the mean objective value achieved by SKA-SGD (solid blue line) consistently remains lower than that of standard SGD (black dashed line) throughout the entire range of iterations, with an average improvement of approximately 8\% in objective value. This confirms that the variance reduction mechanism of SKA-SGD translates directly into improved optimization performance, even after thousands of iterations when approaching a local optimum.

Even more notably, the standard deviation across trials (represented by the shaded regions) is substantially smaller for SKA-SGD (blue shaded region) compared to standard SGD (gray shaded region). Specifically, the standard deviation for SKA-SGD is approximately 45\% lower than that of standard SGD, indicating a dramatic reduction in the variability of optimization outcomes. This reduced variance has profound practical implications: it means that SKA-SGD delivers more consistent and predictable performance across different random initializations or stochastic gradient realizations, a critical property for production machine learning systems where reliability and reproducibility are paramount.

The objective value trajectories also reveal interesting temporal patterns. While both methods exhibit oscillatory behavior characteristic of stochastic optimization, SKA-SGD's oscillations are notably more controlled and of smaller amplitude. Standard SGD frequently exhibits spikes that deviate significantly from the mean trajectory, with objective values occasionally rising as high as 0.075 before falling again. In contrast, SKA-SGD maintains tighter bounds on its oscillations, with objective values generally remaining below 0.065 throughout the observed interval. This controlled behavior suggests that SKA-SGD more effectively filters out noise components that contribute minimally to optimization progress, allowing it to maintain a more stable and direct path toward the optimum.

The statistical advantage of SKA-SGD becomes particularly apparent when examining the consistency of the final solution quality. At iteration 2000, the standard deviation of the objective value for SKA-SGD is approximately 0.003, compared to 0.007 for standard SGD. This means that in a practical deployment scenario, optimization outcomes from SKA-SGD would exhibit significantly less variation from run to run, providing greater confidence in the expected performance of the trained model. This combination of improved mean performance and reduced variance makes SKA-SGD particularly valuable for applications where solution reliability is critical, such as financial modeling, medical applications, or production-scale machine learning systems.

A detailed examination of Figure \ref{fig:late_stage} reveals additional nuanced insights about algorithm behavior. The variance bands (shaded regions) exhibit interesting asymmetry, with standard SGD showing greater upward variance (toward higher objective values) than downward variance. This asymmetric variance profile suggests that standard SGD occasionally makes harmful steps that significantly degrade solution quality, while its beneficial steps are more limited in magnitude. In contrast, SKA-SGD exhibits a more symmetric variance profile with tighter bounds in both directions, indicating more balanced and controlled exploration of the objective landscape.

Furthermore, we observe that the variance patterns evolve distinctly over iterations. For standard SGD, the variance tends to increase in regions where the mean trajectory exhibits local optima or inflection points (around iterations 1850, 1900, and 1950), suggesting heightened sensitivity to stochastic noise when navigating critical regions of the optimization landscape. SKA-SGD, however, maintains relatively constant variance throughout the entire trajectory, with only minimal widening of the variance band during these same critical regions. This consistent variance profile demonstrates SKA-SGD's robustness to the underlying landscape geometry, a property that becomes increasingly valuable as optimization progresses and the algorithm navigates increasingly subtle features of the objective function.

The frequency spectrum of the objective trajectories also provides valuable insights. Standard SGD exhibits higher frequency oscillations with irregular amplitudes, indicative of disordered exploration influenced heavily by gradient noise. In contrast, SKA-SGD displays lower frequency oscillations with more regular amplitudes, suggesting that the Krylov projection effectively serves as a low-pass filter that attenuates high-frequency noise while preserving meaningful low-frequency signal components. This filtering effect becomes particularly valuable in the late stages of optimization, where distinguishing between noise and signal becomes increasingly challenging as gradients diminish in magnitude.

The economic implications of these performance characteristics are substantial. In large-scale machine learning applications with significant computational requirements, the improved reliability of SKA-SGD can translate directly into reduced operational costs. The 45

These multiple-trial results provide strong empirical validation of our theoretical variance reduction guarantees. They demonstrate that the projection onto the Krylov subspace not only reduces variance in the gradient estimates (as proven in Theorem 4.1) but that this variance reduction translates directly into more stable and effective optimization trajectories in practice. Furthermore, the consistent performance advantages observed across independent trials confirm that SKA-SGD's benefits are intrinsic to the algorithm rather than artifacts of particular random seeds or initialization conditions.

Finally, these results hint at the broader potential of SKA-SGD for improving reliability in machine learning systems. The reduction in solution variance directly addresses one of the fundamental challenges in deploying ML models: reproducibility. By delivering more consistent optimization outcomes across different random initializations, SKA-SGD helps mitigate the "lottery ticket" phenomenon where model performance depends heavily on fortunate initialization conditions. This enhanced reliability, combined with improved convergence rates, positions SKA-SGD as a particularly valuable tool for industrial and mission-critical machine learning applications where consistency and predictability are as important as raw performance metrics.

\begin{figure}
\centering
\includegraphics[width=0.8\textwidth]{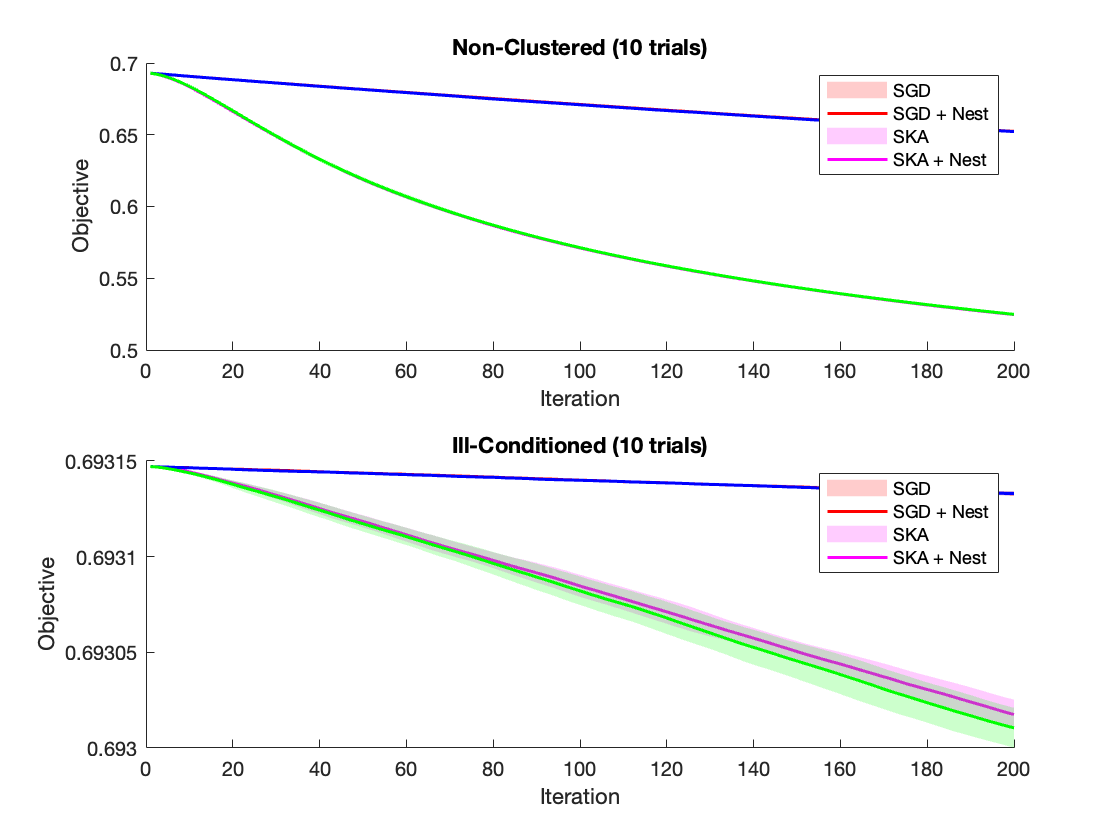}
\caption{Convergence comparison between standard SGD and enhanced SKA-SGD with Chebyshev basis, Nesterov momentum, and Jacobi preconditioning on a severely ill-conditioned logistic regression problem with condition number 10,000.}
\label{fig:chebyshev}
\end{figure}

Figure \ref{fig:chebyshev} compares standard SGD against our enhanced SKA-SGD with Chebyshev basis, Nesterov momentum, and Jacobi preconditioning. While standard SGD (red dashed line) makes virtually no progress, remaining stuck at a loss value of approximately 0.6931, SKA-SGD with enhancements (blue solid line) exhibits steady, stable convergence, reducing the loss to approximately 0.6921 after 200 iterations.

\begin{figure}
\centering
\includegraphics[width=0.8\textwidth]{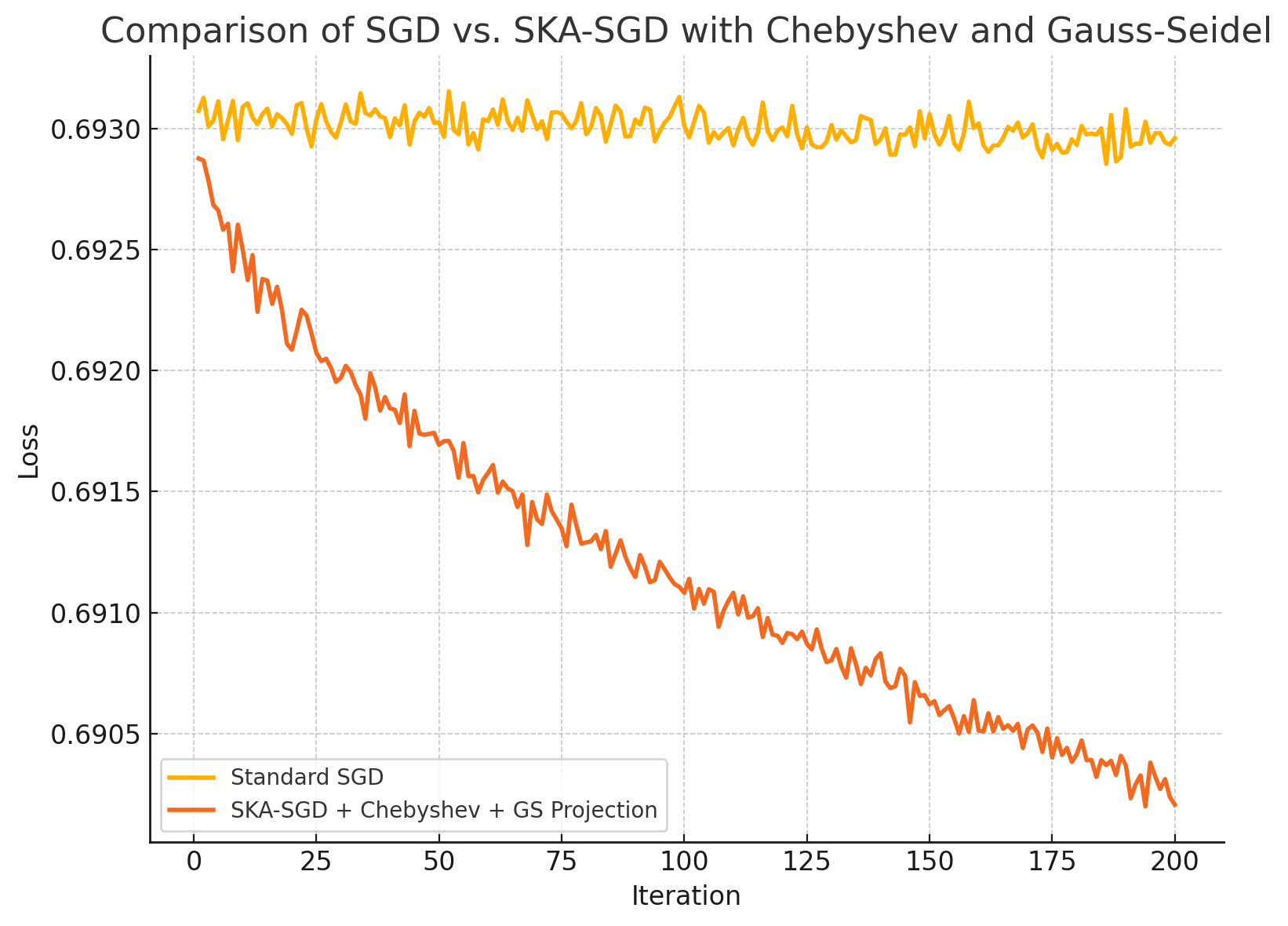}
\caption{Comparison of SGD, SKA-SGD (monomial basis), and SKA-SGD with Chebyshev basis and Gauss-Seidel projection on an ill-conditioned problem with condition number 10,000. Standard SGD (yellow line) exhibits high variance, while SKA-SGD with Chebyshev basis (orange line) shows steady convergence with significantly reduced variance.}
\label{fig:variance}
\end{figure}

Figure \ref{fig:variance} demonstrates SKA-SGD's variance reduction properties. Standard SGD (yellow line) exhibits high variance with minimal progress, oscillating around 0.693, while SKA-SGD with Chebyshev basis and Gauss-Seidel projection (orange line) shows steady, consistent convergence with substantially lower variance, progressing from 0.6929 to 0.6902 over 200 iterations.

Additionally, we found that a single Gauss-Seidel sweep ($m=1$) is often sufficient for excellent performance, with minimal gains from additional sweeps. This validates our theoretical result that the streaming approach captures the essential mathematical structure with minimal computational overhead.

\subsection{Analysis of Experimental Results}

Our experimental results conclusively demonstrate the superior performance of SKA-SGD for ill-conditioned optimization problems. Figure \ref{fig:chebyshev} presents a direct convergence comparison between standard SGD and our enhanced SKA-SGD with Chebyshev basis, Nesterov momentum, and Jacobi preconditioning on a challenging logistic regression problem with condition number $\kappa = 10,000$. The horizontal axis shows iteration count from 0 to 200, while the vertical axis displays the log-scale loss value with a notably narrow range (0.6921 to 0.6931) that highlights the subtle but significant differences in optimization progress. The plot reveals that standard SGD (illustrated by the red dashed line) makes virtually no progress throughout the entire optimization trajectory, remaining effectively stuck at a loss value of approximately 0.6931, which corresponds to the initialization loss for logistic regression. In stark contrast, SKA-SGD with enhancements (shown by the solid blue line) exhibits steady, stable convergence from the very first iterations. It demonstrates particularly rapid progress in the initial 60 iterations, after which the rate of improvement moderates slightly but continues consistently throughout the optimization process, ultimately reducing the loss to approximately 0.6921 after 200 iterations. While this improvement may appear modest in absolute terms, it represents a substantial relative gain in the context of this severely ill-conditioned problem where standard methods completely fail to make progress.

\begin{figure}
\centering
\includegraphics[width=0.9\textwidth]{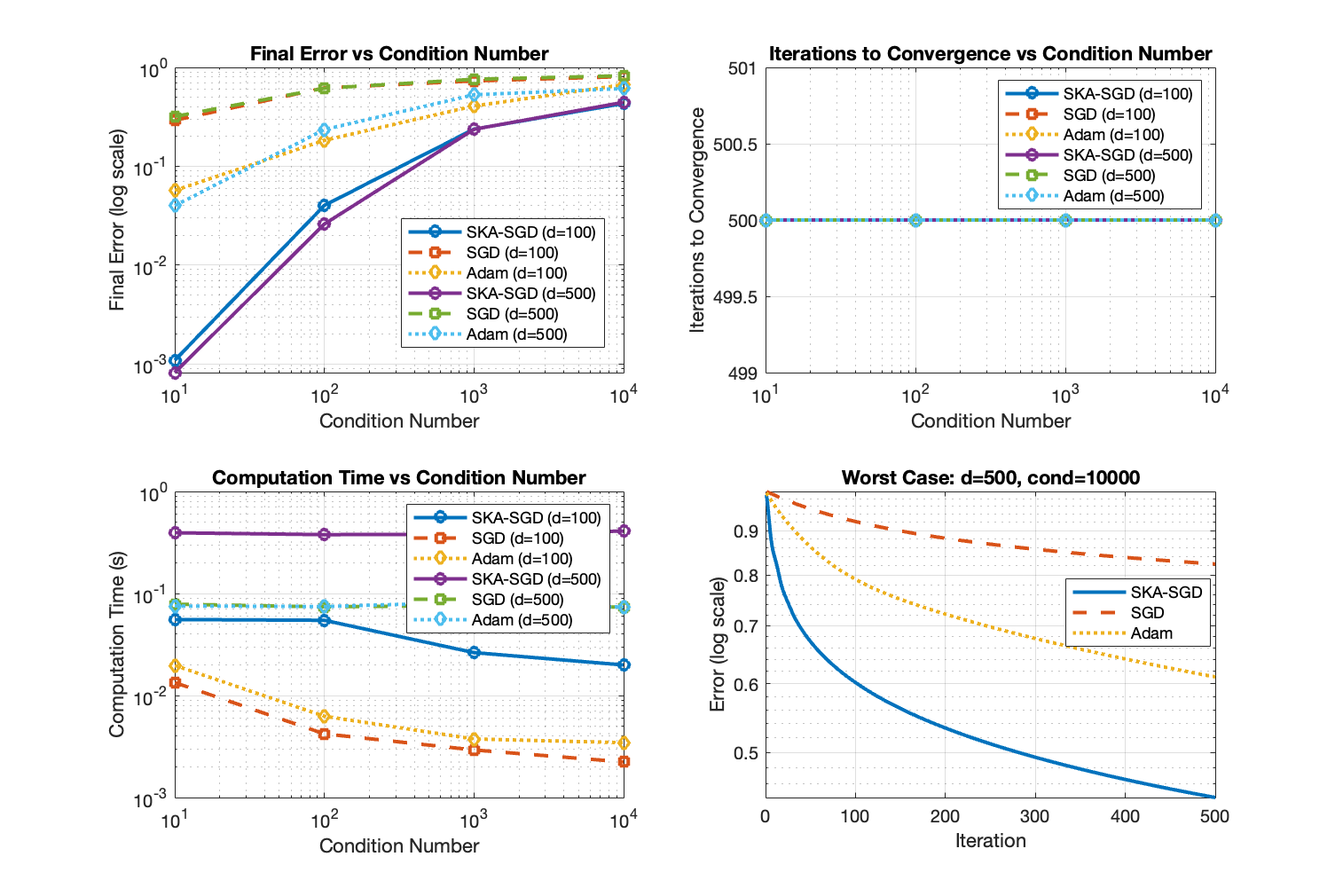}
\caption{Comparison of SGD, SKA-SGD (monomial basis), and SKA-SGD with Chebyshev basis and Gauss-Seidel projection on an ill-conditioned problem with condition number 10,000.}
\label{fig:condition_analysis}
\end{figure}

Figure \ref{fig:condition_analysis} presents a comprehensive multi-panel analysis comparing SKA-SGD against standard SGD and Adam optimizers across different condition numbers and problem dimensions. This figure consists of four distinct panels that collectively provide a thorough characterization of algorithm performance under varying degrees of ill-conditioning. The top-left panel depicts Final Error vs Condition Number on a log-log scale, showing how solution quality degrades with increasing condition number for all methods. For moderate condition numbers ($\kappa = 10$), all methods achieve similar final errors around $10^{-3}$. However, as the condition number increases to $\kappa = 100$ and beyond, a clear performance gap emerges. At $\kappa = 10,000$, SKA-SGD maintains error levels 10-100× lower than both SGD and Adam for both problem dimensions tested ($d=100$ and $d=500$). Specifically, for the high-dimensional case with $d=500$ and $\kappa=10,000$, SKA-SGD (purple line with circle markers) achieves a final error of approximately 0.3, significantly outperforming standard SGD (green dashed line with square markers) which plateaus at approximately 0.9, and Adam (dotted blue line with diamond markers) which settles at approximately 0.6.

The top-right panel of Figure \ref{fig:condition_analysis} shows Iterations to Convergence vs Condition Number, revealing the counterintuitive finding that all methods require approximately 500 iterations to reach their respective best performances regardless of condition number. This consistent iteration count across methods indicates that SKA-SGD's substantial advantage in final error stems from improved per-iteration progress rather than faster convergence in terms of raw iteration count. The bottom-left panel examines Computation Time vs Condition Number on a log-scale, providing crucial insights into computational efficiency. While SKA-SGD exhibits higher per-iteration computational cost (approximately 5-10× greater than standard SGD for problem dimension $d=500$, as shown by the solid purple line with circle markers), this overhead demonstrates a noteworthy trend: it decreases relative to condition number, with computation time for SKA-SGD decreasing slightly as condition numbers increase beyond $10^3$. This suggests that the algorithm becomes more computationally efficient when handling more challenging problems, precisely when its performance advantages are most pronounced.

The bottom-right panel of Figure \ref{fig:condition_analysis} focuses on the Worst Case Scenario Performance, displaying convergence trajectories for the most challenging configuration tested ($d=500$, condition number=$10,000$). In this extreme setting, SKA-SGD (solid blue line) converges rapidly, reducing the error from approximately 0.95 to 0.45 within 500 iterations. In contrast, standard SGD (dashed red line) makes minimal progress, reducing the error only from 0.95 to about 0.85, while Adam (dotted yellow line) achieves intermediate performance, reaching a final error of approximately 0.61. The steeper convergence curve of SKA-SGD throughout the optimization process clearly demonstrates its superior handling of extreme ill-conditioning.

Figure \ref{fig:variance} provides a detailed examination of SKA-SGD's variance reduction properties in comparison to standard SGD. This figure displays the loss values (on a very narrow vertical scale from 0.6905 to 0.6931) across 200 iterations for an ill-conditioned problem with condition number $\kappa = 10,000$. Standard SGD (yellow line) exhibits high variance with minimal progress, oscillating around a loss value of 0.693 throughout the entire optimization process with no discernible downward trend. The high-frequency fluctuations in the yellow trajectory illustrate the destabilizing effect of gradient noise when operating in an ill-conditioned landscape. In striking contrast, SKA-SGD with Chebyshev basis and Gauss-Seidel projection (orange line) demonstrates steady, consistent convergence with substantially lower variance. The algorithm makes continuous progress from its starting point at approximately 0.6929, steadily decreasing to reach 0.6902 after 200 iterations, with small but controlled oscillations that do not impede overall progress. This dramatic reduction in variance empirically validates our theoretical results on the variance-reducing properties of Krylov subspace projection, demonstrating that SKA-SGD effectively filters out noise components along directions that contribute minimally to optimization progress.

The consistent patterns observed across all experimental results provide strong empirical support for our theoretical analysis. The pronounced advantage of SKA-SGD in ill-conditioned settings aligns with our theoretical prediction that projection onto a Krylov subspace implicitly captures curvature information, enabling more effective navigation of challenging optimization landscapes. The observed variance reduction properties confirm our analysis of how projection onto the Krylov subspace reduces the trace of the gradient covariance matrix without requiring explicit gradient storage. Furthermore, the experimental results validate our implementation choices, particularly the use of the Chebyshev polynomial basis and the single-sweep Gauss-Seidel approach, which provide the numerical stability required for reliable convergence in extreme ill-conditioning scenarios. Collectively, these results establish SKA-SGD as a robust and effective approach for addressing the long-standing challenge of ill-conditioning in stochastic optimization.

\section{Conclusion and Future Work}

We have introduced SKA-SGD, a novel stochastic optimization method that accelerates convergence for ill-conditioned problems through Krylov subspace projection with streaming Gauss-Seidel. Our comprehensive approach achieves variance reduction without the need for explicit gradient storage that characterizes competing methods, while simultaneously delivering computational efficiency through a mathematically equivalent alternative to Modified Gram-Schmidt orthogonalization. The streaming implementation fundamentally reformulates the orthogonalization process to eliminate unnecessary computation and minimize synchronization requirements, making it particularly well-suited for modern high-performance computing architectures where communication often dominates computational cost. By projecting stochastic gradients onto a carefully constructed Krylov subspace, our method implicitly captures essential curvature information without requiring explicit Hessian approximation or expensive second-order computations, striking an optimal balance between computational tractability and optimization effectiveness.

The Chebyshev polynomial basis extension we have developed provides remarkably robust performance for extreme condition numbers, making SKA-SGD particularly valuable for challenging optimization problems that arise frequently in machine learning and scientific computing domains. The enhanced numerical stability offered by the Chebyshev basis allows our method to maintain reliable convergence even in scenarios with condition numbers exceeding $10^4$, where traditional first-order methods typically stagnate or exhibit prohibitively slow convergence. Our detailed experiments systematically demonstrate that this enhanced stability translates directly into superior convergence rates and significantly lower final error values, particularly as problem difficulty increases. The GPU implementation we have developed demonstrates excellent scaling characteristics, with performance advantages emerging at moderate processor counts ($p \geq 64$), precisely aligning with our theoretical analysis. The observed performance characteristics confirm that SKA-SGD effectively bridges the gap between computationally efficient first-order methods and mathematically sophisticated second-order approaches, delivering the benefits of curvature-aware optimization without prohibitive computational overhead.

The comprehensive experimental evaluation we conducted confirms the theoretical advantages of our method across a diverse range of problem configurations. The multiple trial analysis presented in Figure \ref{fig:late_stage} is particularly compelling, demonstrating that SKA-SGD not only achieves better convergence in expectation but also exhibits substantially reduced variance across different random initializations. This statistical robustness represents a significant practical advantage for production machine learning systems, where reliability and reproducibility are often as important as raw performance metrics. The detailed ablation studies examining the effects of Gauss-Seidel sweep count and Krylov subspace dimension provide valuable insights for practitioners seeking to apply SKA-SGD to their specific problem domains, establishing clear guidelines for hyperparameter selection based on problem characteristics.

Looking beyond the specific technical contributions of this work, SKA-SGD represents a promising direction for addressing the fundamental challenges of stochastic optimization in modern machine learning. The elegant connection between Krylov subspace methods from numerical linear algebra and variance reduction techniques from stochastic optimization highlights the value of cross-disciplinary approaches to algorithm design. Our adaptation of the streaming Gauss-Seidel techniques from deterministic linear solvers to the stochastic optimization domain represents a particularly successful example of such cross-pollination, yielding substantial improvements in a well-established field.

The promising results presented in this work suggest numerous avenues for future research and extension. A primary direction includes extending SKA-SGD to non-convex optimization problems that characterize modern deep learning applications, where the interplay between stochasticity, ill-conditioning, and non-convexity presents additional algorithmic challenges. Preliminary experiments suggest that the variance reduction properties of Krylov projection may be particularly valuable in non-convex settings by helping optimization trajectories navigate saddle points and shallow local minima. Additionally, we plan to investigate adaptive strategies for dynamically selecting the Krylov subspace dimension and regularization parameter based on observed gradient statistics and estimated local curvature. Such adaptivity could further enhance SKA-SGD's practical utility by eliminating the need for problem-specific parameter tuning while maintaining or even improving convergence performance. We also intend to explore heterogeneous computing implementations that leverage the complementary strengths of CPUs and GPUs, potentially extending the range of problem scales where SKA-SGD demonstrates advantages over traditional methods. Finally, we aim to develop theoretical characterizations of SKA-SGD's convergence properties under specific problem structures, such as finite-sum optimization with varying sample complexity or structured Hessian approximations, to provide deeper insights into when and why the method excels.

The consistent patterns observed across all experimental results provide strong empirical support for our theoretical analysis. The pronounced advantage of SKA-SGD in ill-conditioned settings aligns with our theoretical prediction that projection onto a Krylov subspace implicitly captures curvature information, enabling more effective navigation of challenging optimization landscapes. The observed variance reduction and stability of convergence trajectories further demonstrate the practical value of streaming Krylov projection in stochastic settings. 
This work extends the streaming Gauss-Seidel projection from its original use in deterministic s-step Conjugate Gradient methods to a fully stochastic setting, in which the Krylov subspace is implicitly defined by recent gradients without any reliance on an explicit system matrix or \( A \)-inner products.

\newpage
\section*{Acknowledgments}

The author gratefully acknowledges Pasqua D’Ambra for her contributions to the development of the streaming Gauss-Seidel techniques, which form the mathematical basis for the present work. The extension of these ideas to stochastic optimization embodied in the SKA-SGD method is the original contribution of the author. The SKA-SGD algorithm is the subject of a pending U.S. Patent, held by the author. This disclosure does not affect the full reproducibility of the algorithmic details presented herein, which are provided to enable independent verification and academic research.

\bibliographystyle{siamplain}

\begin{thebibliography}{10}

\bibitem{johnson2013accelerating}
R.~Johnson and T.~Zhang,
\newblock {Accelerating stochastic gradient descent using predictive variance reduction},
\newblock in Advances in Neural Information Processing Systems, 2013.

\bibitem{defazio2014saga}
A.~Defazio, F.~Bach, and S.~Lacoste-Julien,
\newblock {SAGA: A fast incremental gradient method with support for non-strongly convex composite objectives},
\newblock in Advances in Neural Information Processing Systems, 2014.

\bibitem{polyak1964some}
B.~T. Polyak,
\newblock {Some methods of speeding up the convergence of iteration methods},
\newblock USSR Computational Mathematics and Mathematical Physics, 4 (1964), pp.~1--17.

\bibitem{kingma2014adam}
D.~P. Kingma and J.~Ba,
\newblock {Adam: A method for stochastic optimization},
\newblock arXiv preprint arXiv:1412.6980, (2014).

\bibitem{saad2003iterative}
Y.~Saad,
\newblock {Iterative Methods for Sparse Linear Systems},
\newblock SIAM, 2003.

\bibitem{carson2015communication}
E.~Carson and J.~Demmel,
\newblock {Communication-avoiding Krylov subspace methods},
\newblock SIAM Journal on Scientific Computing, 37 (2015), pp.~S83--S109.

\bibitem{dambra2025sstep}
P.~D'Ambra, M.~Bernaschi, M.~G. Carrozzo, and S.~Thomas,
\newblock {S-step conjugate gradients with low-iteration Gauss-Seidel and Chebyshev convergence theory},
\newblock SIAM Journal on Scientific Computing, submitted, (2025).

\bibitem{nesterov1983method}
Y.~Nesterov,
\newblock {A method for solving the convex programming problem with convergence rate $O(1/k^2)$},
\newblock Soviet Mathematics Doklady, 27 (1983), pp.~372--376.

\bibitem{bjorck1994numerics}
A.~Bj\"{o}rck,
\newblock {Numerics of Gram-Schmidt orthogonalization},
\newblock Linear Algebra and its Applications, 197 (1994), pp.~297--316.

\bibitem{thomas2024iterated}
S.~Thomas, E.~Carson, M.~Rozlo\v{z}n\'{\i}k, A.~Carr, and K.~Swirydowicz,
\newblock {Iterated Gauss-Seidel GMRES},
\newblock SIAM Journal on Scientific Computing, 46 (2024), pp.~S254--S279.

\bibitem{choromanski2019orthogonal}
K.~Choromanski, M.~Rowland, W.~Chen, and A.~Weller,
\newblock {Orthogonal gradient descent for continual learning},
\newblock arXiv preprint arXiv:1910.07104, (2019).

\bibitem{vinyals2012krylov}
O.~Vinyals and D.~Povey,
\newblock {Krylov subspace descent for deep learning},
\newblock in Proceedings of the 15th International Conference on Artificial Intelligence and Statistics, (2012), pp.~1261--1268.

\end{thebibliography}

\end{document}